\documentclass[10pt]{article}
\usepackage{amsmath}
\usepackage{amssymb}
\usepackage{amsthm}
\usepackage{mathrsfs}
\usepackage{tikz,authblk}
\newtheorem{theorem}{Theorem}[section]
\newtheorem{lemma}[theorem]{Lemma}
\newtheorem{corollary}[theorem]{Corollary}
\newtheorem{proposition}[theorem]{Proposition}
\newtheorem{conjecture}[theorem]{Conjecture}

\theoremstyle{definition}

\begin{document}
\title{The second out-neighbourhood for local tournaments\thanks{Corresponding author. $E$-$mail\ \ address:$ ruijuanli@sxu.edu.cn(R. Li). Research of RL is partially supported by NNSFC under no. 11401353 and TYAL of Shanxi.}}
\author{Ruijuan Li}
\author{Juanjuan Liang}
\affil{School of Mathematical Sciences, Shanxi University, Taiyuan, Shanxi, 030006, PR China}

\maketitle

\begin{abstract}
Sullivan stated the conjectures: (1) every oriented graph $D$ has a vertex $x$ such that
$d^{++}(x)\geq d^{-}(x)$; (2) every oriented graph $D$ has a vertex $x$ such that $d^{++}(x)+d^{+}(x)\geq 2d^{-}(x)$.
In this paper, we prove that these conjectures hold for local tournaments. In particular, for a local tournament $D$, we prove that $D$ has at least two vertices satisfying $(1)$ if $D$ has no vertex of in-degree zero. And, for a local tournament $D$, we prove that  either there exist two vertices satisfying $(2)$ or there exists a vertex $v$ satisfying  $d^{++}(v)+d^{+}(v)\geq 2d^{-}(v)+2$ if $D$ has no vertex of in-degree zero.
\end{abstract}

\vskip 3mm
\noindent{\bf Keywords:} local tournaments; the second out-neighbourhood; Sullivan's conjectures; round decomposable; non-round decomposable
\vskip 3mm



\section{Introduction}

In this paper, we consider finite digraphs without loops and multiple arcs. The main source for terminology and notation is \cite{Bang-Jensen}.

Let $D$ be a digraph. We denote the vertex set and the arc set of $D$ by $V(D)$ and $A(D)$, respectively. For a vertex subset $X$, we denote  the subdigraph of $D$ induced by $X$ (respectively, $D-X$) by $D\langle X\rangle$ (respectively, $D\langle V(D)-X\rangle$). For convenience, we write $D-X$ instead of $D\langle V(D)-X\rangle$. In addition $D-x=D-\{x\}$. And if $X$ is a subdigraph, we write $D-X$ instead of $D-V(X)$.

Let $x,y$ be distinct vertices of $D$. If there is an arc from $x$ to $y$, we say that $x$ {\it dominates} $y$ and denote it by $x \rightarrow y$ and call $y$ (respectively, $x$) an {\it out-neighbour} (respectively, an {\it in-neighbour)} of $x$ (respectively, $y$). If $V_{1}$ and $V_{2}$ are disjoint subsets of vertices of $D$ such that there is no arc from $V_{2}$ to $V_{1}$ and $a\rightarrow b$ for all $a\in V_{1}$ and $b\in V_{2}$, then we say that $V_{1}$ {\it completely dominates} $V_{2}$ and denote it by $V_{1}\Rightarrow V_{2}$. We will use the same notation when $V_{1}$ or $V_{2}$ is subdigraphs of $D$. In particular, if $V_{1}$ contains only one vertex $v$,  denote it by $v\Rightarrow V_{2}$.

For a subdigraph or simply a vertex subset $H$ of $D$ (possibly, $H = D$), we let $N^{+}_{H}(x)$ (respectively, $N^{-}_{H}(x)$) denote the set of out-neighbours (respectively, in-neighbours) of $x$ in $H$ and call it the {\it out-neighbourhood} (respectively, {\it in-neighbourhood}) of $x$ in $H$. Furthermore, $d^{+}_{H}(x) = |N^{+}_{H}(x)|$ (respectively, $d^{-}_{H}(x) = |N^{-}_{H}(x)|)$ is called the {\it out-degree} (respectively, {\it in-degree}) of $x$ in $H$. Let$$N^{++}_{H}(x) =\bigcup\limits_{u\in N^{+}_{H}(x)} N^{+}_{H}(u)-N^{+}_{H}(x)$$ which is called the {\it second out-neighbourhood} of $x$ in $H$. Furthermore, $d^{++}_{H}(x) = |N^{++}_{H}(x)|$. We will omit the subscript $H$ if the digraph is known from the context. For a pair of vertex disjoint subdigraphs $H$ and $H'$, we define $$N_{H}^{+}(H')=\bigcup\limits_{x\in V(H')} N^{+}_{H}(x)-V(H'),\,\,\,\, N_{H}^{-}(H')=\bigcup\limits_{x\in V(H')} N^{-}_{H}(x)-V(H').$$

A vertex $x$ is a {\it 2-king} (for short, a {\it king}) of $D$, if for $y\in V(D)-x$, there exists an $(x,y)$-path of length at most 2.

A digraph $D$ is {\it strong} if, for every pair $x$, $y$ of distinct vertices, $D$ contains a  path from $x$ to $y$ and a
path from $y$ to $x$. A strong component of a digraph $D$ is a maximal induced subdigraph of $D$ which is strong. If $D_{1},D_{2},\ldots,D_{t}$ are the {\it strong components} of $D$, then clearly $V(D_{1})\cup V(D_{2})\cup \ldots \cup V(D_{t}) =V(D)$ (note that a digraph with only one vertex is strong). Moreover, we must have $V(D_{i})\cap V(D_{j})=\emptyset$ for every $i\neq j$. The strong components of $D$ can be labelled $D_{1},D_{2},\ldots,D_{t}$ such that there is no arc from $D_{j}$ to $D_{i}$ unless $j<i$. We call such an ordering an {\it acyclic ordering of the strong components} of $D$.

For a vertex subset $S$ of strong digraph $D$, $S$ is called a {\it separating set} of $D$ if $D-S$ is not strong. A separating set $S$ of $D$ is {\it minimal} if for any proper subset $S'$ of $S$, the subdigraph $D-S'$ is strong.

A digraph  $R$ on $n$ vertices is {\it round} if we can label its vertices $v_{1},v_{2},\ldots,v_{n}$ so that for each $i$, we have
$N^{+}(v_{i})=\{v_{i+1},v_{i+2},\ldots,v_{i+d^{+}(v_{i})}\}$ and $N^{-}(v_{i})=\{v_{i-d^{-}(v_{i})},\ldots,v_{i-2},v_{i-1}\}$(all subscripts are taken modulo $n$). We will refer to the ordering $v_{1},v_{2},\ldots,v_{n}$ as a {\it round labelling} of $R$.

A digraph $D$ is {\it semicomplete} if, for every pair $x$, $y$ of distinct vertices in $D$, either $x$ dominates $y$ or $y$ dominates $x$ (or both). {\it Tournaments} are semicomplete digraphs with no 2-cycle.

A digraph $D$ with no 2-cycle is an {\it oriented graph}.

In 1990, Seymour \cite{Dean} proposed the following conjecture which is one of the most interesting and challenging open questions concerning oriented graphs.

\begin{conjecture}\label{1.1}(Seymour's Second Neighbourhood Conjecture(SSNC))
For any oriented graph $D$, there exists a vertex $v$ in $D$ such that $d^{++}(v)\geq d^{+}(v)$.
\end{conjecture}

We call such a vertex $v$ satisfying Conjecture \ref{1.1} a {\it Seymour vertex}. The first non-trivial result for SSNC was obtained by Fisher \cite{Fisher} who proved Dean's conjecture \cite{Dean}, which is SSNC restricted to tournaments. Fisher used Farkas' Lemma and averaging arguments.

\begin{theorem}\label{1.2} \cite{Fisher} In any tournament $T$, there exists a Seymour vertex.
\end{theorem}

A more elementary proof of SSNC for tournaments was given by Havet and
Thomass\'{e} \cite{Havet} who introduced a median order approach. Their proof also yields
the following stronger result.

\begin{theorem}\label{1.3}\cite{Havet} A tournament $T$ with no vertex of out-degree zero has at
least two Seymour vertices.
\end{theorem}

Fidler and Yuster \cite{Fidler2} further developed the median order approach and proved that SSNC holds for oriented graphs $D$ with minimum degree $|V(D)|-2$, tournaments minus a star and tournaments minus the arc set of a subtournament. The median order approach was also used by Ghazal \cite{Ghazal} who proved a weighted version of SSNC for tournaments missing a generalized star. Kaneko and Locke \cite{Kaneko} proved SSNC for oriented graphs with minimum out-degree at most $6$. Cohn, Godbole, Wright, Harkness and Zhang \cite{Cohn} proved SSNC for random oriented graphs with probability $p<\frac{1}{2}-\delta$. Gutin and Li \cite{Gutin} proved SSNC for extended tournaments and quasi-transitive oriented graphs.

Another approach to SSNC is to determine the maximum value $\gamma$ such that in every oriented graph $D$, there exists a vertex $x$ such that $d^{+}(x)\leq \gamma d^{++}(x)$. SSNC asserts that $\gamma=1$. Chen, Shen and Yuster \cite{Chen} proved that $\gamma \geq r$ where $r=0.657298\ldots$ is the unique real root of $2x^{3}+x^{2}-1=0$. Furthermore, they improves this bound to $0.67815\ldots$  mentioned in the end of the article \cite{Chen}.

\vskip 2mm

Sullivan \cite{Sullivan} stated the following ``compromise conjectures'' on SSNC, where $d^{-}(v)$ is used instead of or together with $d^{+}(v)$.

\begin{conjecture}\label{1.4} \cite{Sullivan}
 (1) Every oriented graph $D$ has a vertex $x$ such that $d^{++}(x)\geq d^{-}(x)$.

 (2) Every oriented graph $D$ has a vertex $x$ such that $d^{++}(x)+d^{+}(x)\geq 2d^{-}(x)$.
\end{conjecture}

For convenience, a vertex $x$ in $D$ satisfying Conjecture \ref{1.4} ($i$) is called a {\it Sullivan-$i$} vertex of $D$ for $i\in\{1,2\}$.

Li and Sheng \cite{Li quasi} \cite{Li tournament}, proved Sullivan's Conjectures for tournaments, extended tournaments, quasi-transitive oriented graphs as well as bipartite tournaments. For tournaments, they obtained the following results:

\begin{corollary}\label{1.5}\cite{Li quasi} Every tournament has a Sullivan-$1$ vertex and a Sullivan-$2$ vertex. Every tournament with no vertex of in-degree zero has at least three Sullivan-1 vertices.
\end{corollary}

\begin{theorem}\label{1.61}\cite{Li quasi} A tournament $T$ has at least two Sullivan-$2$ vertices unless $T\in \mathcal{T}$.
\end{theorem}

$\mathcal{T}$  is a special class of tournaments. $T\in \mathcal{T}$ if $T$ is a tournament consisting of exactly two strong components $T_{1}$ and $T_{2}$ such that $T_{1}$ dominates $T_{2}$, $T_{1}$ is a single vertex $v$ and $T_{2}$ is a tournament satisfying that $d_{T_{2}}^{+}(x)\leq d_{T_{2}}^{-}(x)+1$ for any $x\in V(T_{2})$. It is easy to check that $v$ is the unique Sullivan-2 vertex of $T$.

\vskip 2mm

From Theorem \ref{1.61}, we obtain immediately the following result:

\begin{corollary}\label{1.6}
A strong tournament $T$ with at least three vertices has at least two Sullivan-$2$ vertices.
\end{corollary}

A digraph $D$ is {\it locally semicomplete}  if $D\langle N^{+}(x)\rangle$ and $D\langle N^{-}(x)\rangle$ are both semicomplete for every vertex $x$ of $D$. Specifically, every round digraph is locally semicomplete \cite{Guo locally}. A {\it local tournament} is a locally semicomplete digraph with no 2-cycle.

Let $D$ be a digraph with vertex set $\{v_{1},v_{2},\ldots,v_{n}\}$, and let $G_{1},G_{2},\ldots,G_{n}$ be
digraphs which are pairwise vertex disjoint. The {\it composition} $D[G_{1},G_{2},\ldots,G_{n} ]$
is the digraph $L$ with vertex set $V(G_{1})\cup V(G_{2})\cup\ldots\cup V(G_{n})$ and arc set
$(\cup^{n}_{i=1}A(G_{i}))\cup \{g_{i}g_{j}|g_{i}\in V(G_{i}), g_{j}\in V(G_{j}), v_{i}v_{j}\in A(D)\}$. If $D=H[V_{1},V_{2},\ldots,\\V_{n}]$ and none of the digraphs $V_{1},V_{2},\ldots,V_{n}$ has an arc, then $D$ is an {\it extension} of $H$.

A digraph $D$ is {\it round decomposable} if there exists a round local tournament $R$ on $r\geq 2$
vertices such that $D=R[S_{1},S_{2},\ldots,S_{r}]$, where each $S_{i}$ is a strong semicomplete digraph. We call $R[S_{1},S_{2},\ldots,S_{r}]$ a {\it round decomposition} of $D$. Clearly, a round decomposable digraph is locally semicomplete.

Locally semicomplete digraphs were introduced in 1990 by Bang-Jensen \cite{Bj}. The following theorem, due to Bang-Jensen,
Guo, Gutin and Volkmann, stated a full classification of locally semicomplete digraphs.

\begin{theorem}\label{1.7}\cite{Guo locally} Let $D$ be a connected locally semicomplete digraph. Then exactly one of the following possibilities holds:

$(a)$ $D$ is round decomposable with a unique round decomposition $R[S_{1},S_{2},\ldots,\\S_{r}]$, where $R$ is a round local tournament on $r\geq2$ vertices and $S_{i}$ is a strong semicomplete digraph for each $i\in \{1,2,\ldots,r\}$;

$(b)$ $D$ is non-round decomposable and not semicomplete and it has the structure as described in Theorem \ref{2.5};

$(c)$ $D$ is a semicomplete digraph which is non-round decomposable.
\end{theorem}

If $D$ is restricted to a local tournament, we have the following result:

\begin{corollary}\label{1.8} Let $D$ be a connected local tournament. Then exactly one of the following possibilities holds:

$(a)$ $D$ is round decomposable with a unique round decomposition $R[S_{1},S_{2},\ldots,\\ S_{r}]$, where $R$ is a round local tournament on $r\geq2$ vertices and $S_{i}$ is a strong tournament for $i\in \{1,2,\ldots,r\}$;

$(b)$ $D$ is non-round decomposable and not a tournament and it has the structure as described in Theorem \ref{2.5};

$(c)$ $D$ is a tournament which is non-round decomposable.
\end{corollary}

In \cite{Li seymour}, we investigate SSNC for local tournaments. In this paper, we discuss Sullivan's Conjectures for local tournaments. In Section 2, we introduce the structure of a local tournament. In Section 3 and Section 4, we investigate the Sullivan-$i$ vertex in a round decomposable local tournament and a non-round decomposable local tournament, respectively, for $i\in\{1,2\}$.


\section{The structure of a local tournament}

In this section, all theorems  are on the structure of locally semicomplete digraphs. Clearly,
these theorems also hold if the digraph is restricted to a local tournament.

\vskip 1mm

\begin{theorem}\label{2.1}\cite{Con} Let $D$ be a connected, but not strong locally semicomplete digraph. Then the following holds for $D$.

$(a)$ If $A$ and $B$ are distinct strong components of $D$ with at least one arc between them, then either $A\Rightarrow B$ or $B\Rightarrow A$.

$(b)$ If $A$ and $B$ are strong components of $D$, such that $A\Rightarrow B$, then $A$ and $B$ are semicomplete digraphs.

$(c)$ The strong components of $D$ can be ordered in a unique way $D_{1},D_{2},\ldots,D_{p}$ such that there is no arc from $D_{j}$ to $D_{i}$ for $j > i$, and $D_{i}$ completely dominates $D_{i+1}$ for $i\in \{1,2,\ldots,p-1\}$.
\end{theorem}

A kind of the decomposition of non-strong locally semicomplete digraphs described in \cite{Con} is the following.

\begin{theorem}\label{2.2}\cite{Con} Let $D$ be a connected, but not strong locally semicomplete digraph, and let $D_{1},D_{2},\ldots,D_{p}$ be the acyclic ordering of the strong components of $D$. Then $D$ can be decomposed into $r\geq2$ induced subdigraphs $D'_{1},D'_{2},\ldots,D'_{r}$ as follows:
$$D'_{1}=D_{p}, \,\,\,\lambda_{1}=p\mbox{,}$$
$$\lambda_{i+1}=\min\{j\,|\,N^{+}(D_{j})\cap V(D'_{i})\neq\emptyset\}\mbox{\,\,for each} \,\,i\in \{1,2,\ldots,r-1\}\mbox{,}$$
$$\mbox{and}\,\,D'_{i+1}=D\langle V(D_{\lambda_{i+1}})\cup V(D_{\lambda_{i+1}+1})\cup \ldots \cup V(D_{\lambda_{i}-1})\rangle \mbox{.}$$

The subdigraphs $D'_{1},D'_{2},\ldots,D'_{r}$ satisfy the properties below:

$(a)$ $D'_{i}$ consists of some strong components of $D$ and is semicomplete for each $i\in \{1,2,\ldots,r\}$;

$(b)$ $D'_{i+1}$ completely dominates the initial component of $D'_{i}$ and there exists no arc from $D'_{i}$ to $D'_{i+1}$ for $i \in \{1,2,\ldots,r-1\}$;

$(c)$ If $r\geq 3$, then there is no arc between $D'_{i}$ and $D'_{j}$ for $i,\,j$ satisfying $|j-i|\geq2$.
\end{theorem}

The unique sequence $D'_{1},D'_{2},\ldots,D'_{r}$ defined in Theorem \ref{2.2} will be referred to as the $semicomplete$ $decomposition$ of $D$.

\begin{theorem}\label{2.11}\cite{Guo locally}
If $D$ is a round decomposable locally semicompete digraph, then it has a unique round decomposition $D=R[S_{1},S_{2},\ldots,S_{r}]$, where $R$ is a round local tournament on $r\geq 2$
vertices and each $S_{i}$ is a strong semicomplete digraph.
\end{theorem}

\begin{theorem}\label{2.5}\cite{Guo locally} 
Let $D$ be a strong locally semicomplete digraph which is not semicomplete. Then $D$ is non-round decomposable if and only if the following conditions are satisfied:

$(a)$ There is a minimal separating set $S$ such that $D-S$ is not semicomplete, and for each such $S$, $D\langle S\rangle$ is semicomplete and the
semicomplete decomposition of $D-S$ has exactly three components $D'_{1}$, $D'_{2}$, $D'_{3}$;

$(b)$ There are integers $\alpha, \beta, \mu, \nu$ with $\lambda \leq \alpha \leq \beta \leq p-1$ and $p+1 \leq \mu \leq \nu\leq p+q$ such that
$$N^{-}(D_{\alpha})\cap V(D_{\mu})\neq \emptyset \quad and \quad N^{+}(D_{\alpha})\cap V(D_\nu)\neq \emptyset\mbox{,}$$
$$\mbox{ or }\quad N^{-}(D_\mu)\cap V(D_\alpha)\neq  \emptyset \quad and \quad N^{+}(D_\mu)\cap V(D_\beta)\neq \emptyset\mbox{,}$$
where $D_{1},D_{2},\ldots,D_{p}$ and $D_{p+1},D_{p+2},\ldots,D_{p+q}$ are the acyclic orderings of the strong components of $D-S$ and $D\langle S\rangle$, respectively, and $D_{\lambda}$ is the initial component of $D_{2}^{'}$.\end{theorem}

By Theorem \ref{2.5}, $D$ is always strong if $D$ is a non-round decomposable locally semicomplete digraph. An example of a non-round decomposable locally semicomplete digraph is shown in Figure 2.

\vskip 1mm

\begin{theorem}\label{2.6} \cite{Guo locally} 
Let $D$ be a strong non-round decomposable locally semicomplete digraph and let S be a minimal separating set of
$D$ such that $D-S$ is not semicomplete. Let $D_{1},D_{2},\ldots,D_{p}$ be the acyclic ordering of the strong components of $D-S$ and $D_{p+1}, D_{p+2},\ldots,D_{p+q}$ be the acyclic ordering of the strong components of $D\langle S\rangle$. The following holds:

$(a)$ $D_{p}\Rightarrow S\Rightarrow D_{1}$.

$(b)$ Suppose that there is an arc $s\rightarrow v$ from $S$ to $D'_{2}$ with $s \in V(D_{i})$ and $v\in V(D_{j})$. Then
$D_{i}\cup D_{i+1} \cup \ldots \cup D_{p+q}\Rightarrow D'_{3}\Rightarrow D_{\lambda} \cup\ldots\cup D_{j}$.

$(c)$ $ D_{p+q}\Rightarrow D'_{3}$ and $D_{f}\Rightarrow D_{f+1}$ for $f\in\{1,2,\ldots,p+q\}$ where subscripts are modulo $p + q$.
\end{theorem}

\section{In a round decomposable local tournament}

In this section, $D$ is always a round decomposable local tournament and let the unique round decomposition of $D$ be $R[S_{1},S_{2},\ldots,S_{r}]$,  where $R$ is a round local tournament on $r\geq 2$ vertices and each $S_{i}$ is a strong tournament.

\vskip 3mm

We begin with a useful observation.

\begin{lemma}\label{121} 
Let $D$ be a round decomposable local tournament and $R[S_{1},S_{2},\ldots,\\S_{r}]$ be the unique round decomposition of $D$. Let $D^{*}=R[V_{1},V_{2},\ldots,V_{r}]$ and $v_{i}\in V_{i}$ be arbitrary, where $V_{i}$ is the vertex set of $S_{i}$ for $i\in \{1,2,\ldots,r\}$. If $N_{D^{*}}^{+}(v_{j})=V_{j+1}\cup \ldots\cup V_{k}$, then $d_{D^{*}}^{++}(v_{j})\geq d_{D^{*}}^{+}(v_{k})$.
\end{lemma}
\noindent{\it Proof.} Let $v\in N_{D^{*}}^{+}(v_{k})$. We claim that $v\notin N_{D^{*}}^{+}(v_{j})$. In fact, if $v_{j}\rightarrow v$, then $v_{j}$, $v$, $v_{k}$ are in the order of the round labelling of $R$. Then $v_{k}\rightarrow v_{j}$ since $v_{k}\rightarrow v$. Note that $v_{j}\rightarrow v_{k}$. This contradicts the fact that $D$ has no 2-cycle. So $v\notin N_{D^{*}}^{+}(v_{j})$. Thus $v\in N_{D^{*}}^{++}(v_{j})$, i.e., $N_{D^{*}}^{+}(v_{k})\subseteq N_{D^{*}}^{++}(v_{j})$. Then $d^{++}_{D^{*}}(v_{j})\geq d^{+}_{D^{*}}(v_{k})$.\qed

\vskip 3mm
First, we consider the existence of a Sullivan-$i$ vertex in $D$ for $i\in\{1,2\}$.

\begin{lemma}\label{3.1} 
Let $D$ be a round decomposable local tournament and $R[S_{1},S_{2},\ldots,\\S_{r}]$ be the unique round decomposition of $D$. Let $D^{*}=R[V_{1},V_{2},\ldots,V_{r}]$ and $v_{j}\in V_{j}$ be arbitrary, where $V_{j}$ is the vertex set of $S_{j}$ for $j\in \{1,2,\ldots,r\}$. If there is a vertex $v\in V_{j}$ such that $v$ is a Sullivan-i vertex of $S_{j}$ and a Sullivan-i vertex of $D^{*}$, then $v$ is a Sullivan-i vertex of $D$  for $i\in\{1,2\}$.
\end{lemma}
\noindent{\it Proof.} For the case when $i=1$, since $v$ is a Sullivan-1 vertex of $S_{j}$ and a Sullivan-1 vertex of $D^{*}$, we have $d^{++}_{S_{j}}(v)\geq d^{-}_{S_{j}}(v),\,\,d^{++}_{D^{*}}(v)\geq d^{-}_{D^{*}}(v)$. Clearly $$d^{++}_{D}(v)=d^{++}_{S_{j}}(v)+d^{++}_{D^{*}}(v),\,\,d^{-}_{D}(v)=d^{-}_{S_{j}}(v)\cup d^{-}_{D^{*}}(v)\mbox{.}$$ Thus $d^{++}_{D}(v)\geq d^{-}_{D}(v)$ and $v$ is a Sullivan-1 vertex of $D$.

For the case when $i=2$, it can be proved similarly.\qed

\vskip 3mm

\begin{theorem} \label{4.2} 
Let $D$ be a round decomposable local tournament. Then $D$ has a Sullivan-i vertex for $i\in\{1,2$\}.
\end{theorem}

\noindent{\it Proof.} Let $R[S_{1},S_{2},\ldots,S_{r}]$ be the unique round decomposition of $D$. Let $D^{*}=R[V_{1},V_{2},\ldots,V_{r}]$ and $v_{j}\in V_{j}$ be arbitrary, where $V_{j}$ is the vertex set of $S_{j}$ for $j\in \{1,2,\ldots,r\}$. W.l.o.g., assume that $v_{1}\in V_{1}$ is a vertex of $D^{*}$ with minimum out-degree, i.e., $d_{D^{*}}^{+}(v_{1})=\delta^{+}(D^{*})$. Let $N_{D^{*}}^{+}(v_{1})=V_{2} \cup\ldots\cup V_{t}$.
Since $v_{1}\nrightarrow v_{t+1}$, we have $N_{D^{*}}^{-}(v_{t+1})\subseteq V_{2} \cup\ldots\cup V_{t}= N_{D^{*}}^{+}(v_{1})$. Then
$$d_{D^{*}}^{-}(v_{t+1})\leq d_{D^{*}}^{+}(v_{1})=\delta^{+}(D^{*})\mbox{.}$$
Let $N_{D^{*}}^{+}(v_{t+1})=V_{t+2} \cup \ldots \cup V_{h}$. By Lemma \ref{121}, we have
$$d_{D^{*}}^{++}(v_{t+1})\geq d_{D^{*}}^{+}(v_{h})\geq \delta^{+}(D^{*})\mbox{.}$$
Thus $d_{D^{*}}^{++}(v_{t+1})\geq\delta^{+}(D^{*})\geq d_{D^{*}}^{-}(v_{t+1})$ and $v_{t+1}$ is a Sullivan-1 vertex of $D^{*}$.
Note that $d_{D^{*}}^{+}(v_{t+1})\geq \delta^{+}(D^{*})$. Then $d_{D^{*}}^{++}(v_{t+1})+d_{D^{*}}^{+}(v_{t+1})\geq 2\delta^{+}(D^{*})\geq2d_{D^{*}}^{-}(v_{t+1})$ and $v_{t+1}$ is also a Sullivan-2 vertex of $D^{*}$.

Clearly, all the vertices of $V_{t+1}$ are Sullivan-$i$ vertices of $D^{*}$ for $i\in\{1,2\}$.
By  Corollary \ref {1.5}, the tournament $S_{t+1}$ always has a Sullivan-1 vertex and a Sullivan-2 vertex, say $v_{t+1}$ and $v'_{t+1}$, respectively. By Lemma \ref{3.1}, $v_{t+1}$ is a Sullivan-1 vertex of $D$ and $v'_{t+1}$ is a Sullivan-2 vertex of $D$.\qed

\vskip 3mm

Next, we consider the number of Sullivan-$i$ vertices in a connected round decomposable local tournament  with no vertex of in-degree zero for $i\in\{1,2\}$. Note that every non-strong local tournament is round decomposable. We consider two cases: $(1)$ a connected, but not strong local tournament with no vertex of in-degree zero; $(2)$ a strong round decomposable local tournament.

\vskip 3mm

\begin{theorem}\label{4.3} 
Let $D$ be a connected, but not strong local tournament with no vertex of in-degree zero. Then $D$ has at least three Sullivan-1 vertices and two Sullivan-2 vertices.
\end{theorem}

\noindent{\it Proof.} Let $D_{1},D_{2},\ldots,D_{p}$ be the acyclic ordering of the strong components of $D$. Let $v$ be a Sullivan-1 vertex of $D_{1}$ and $v'$ be a Sullivan-2 vertex of $D_{1}$, i.e., $$d^{++}_{D_{1}}(v)\geq d^{-}_{D_{1}}(v),\,\,d^{++}_{D_{1}}(v')+d^{+}_{D_{1}}(v')\geq 2d^{-}_{D_{1}}(v')\mbox{.}$$
Clearly$$d^{++}_{D}(v)\geq d^{++}_{D_{1}}(v), \,\, d^{-}_{D}(v)=d^{-}_{D_{1}}(v)\mbox{.}$$
$$d^{++}_{D}(v')\geq d^{++}_{D_{1}}(v'), \,\,d^{+}_{D}(v')\geq d^{+}_{D_{1}}(v'),\,\, d^{-}_{D}(v)=d^{-}_{D_{1}}(v')\mbox{.}$$
Thus $d^{++}_{D}(v)\geq d^{-}_{D}(v)$ and $d^{++}_{D}(v')+d^{++}_{D}(v')\geq 2d^{-}_{D}(v')$, i.e., $v$ is a Sullivan-1 vertex of $D$ and $v'$ is a Sullivan-2 vertex of $D$. So the Sullivan-$i$ vertex of $D_{1}$ is always the Sullivan-$i$ vertex of $D$ for $i\in\{1,2\}$.

Since $D$ has no vertex of in-degree zero, we see that $D_{1}$ has no vertex of in-degree zero and $D_{1}$ has at least three vertices. By Corollary \ref{1.5} and \ref{1.6}, $D_{1}$ has at least three Sullivan-1 vertices and two Sullivan-2 vertices. Then these three vertices (respectively, two vertices) are Sullivan-1 (respectively, Sullivan-2) vertices  of $D$.\qed

\begin{theorem}\label{4.4} Let $D$ be a strong round decomposable local  tournament. Then $D$ has at least two Sullivan-i vertices for $i\in\{1,2\}$.
\end{theorem}
\noindent{\it Proof.} Let $D=R[S_{1},S_{2},\ldots,S_{r}]$ be the unique round decomposition.  Let $D^{*}=R[V_{1},V_{2},\ldots,V_{r}]$ and $v_{j}\in V_{j}$ be arbitrary, where $V_{j}$ is the vertex set of $S_{j}$ for $j\in \{1,2,\ldots,r\}$. W.l.o.g., assume that $v_{1}\in V_{1}$ is a vertex of $D^{*}$ with minimum out-degree, i.e., $d_{D^{*}}^{+}(v_{1})=\delta^{+}(D^{*})$.  Let $N_{D^{*}}^{+}(v_{1})=V_{2}\cup\ldots\cup V_{t}$. According to the proof of Theorem \ref {4.2},  a Sullivan-$i$ vertex of $S_{t+1}$ is a Sullivan-$i$ vertex of $D$ for $i\in\{1,2\}$.

For the case when $|V_{t+1}|\geq 2$, by Corollary \ref{1.5} and \ref{1.6}, the strong tournament $S_{t+1}$ has at least three Sullivan-1 vertices and two Sullivan-2 vertices. Then these three vertices (respectively, two vertices) are the Sullivan-1 (respectively, Sullivan-2) vertices of $D$.

For the case when $|V_{t+1}|= 1$ and there exists $v_{h}\notin V_{1}$ such that $d_{D^{*}}^{+}(v_{h})=\delta^{+}(D^{*})$, we can repeat the proof of Theorem \ref{4.2} and obtain a different $``v_{t+1}$". Now the so-called $``v_{t+1}$" is another Sullivan-$i$ vertex for $i\in\{1,2\}$.

Now we consider the case when $|V_{t+1}|= 1$ and there exists  no $v_{h}\notin V_{1}$ such that $d_{D^{*}}^{+}(v_{h})=\delta^{+}(D^{*})$. Then $d^{+}_{D^{*}}(v_{j})>\delta^{+}(D^{*})\geq1$ for any $j\neq1$. According to the proof of Theorem \ref {4.2}, the only vertex $v_{t+1}$ of $V_{t+1}$ is a Sullivan-$i$ vertex of $D$ for $i\in\{1,2\}$. It is sufficient to find another Sullivan-$i$ vertex of $D$ for $i\in\{1,2\}$.

We claim that $v_{2}\rightarrow v_{t+2}$ and $v_{t+2}\notin V_{1}$. In fact, we have $v_{t+2}$ must be in the set $N_{D^{*}}^{+}(v_{2})$ since $d^{+}_{D^{*}}(v_{2})>d^{+}_{D^{*}}(v_{1})$, $|V_{2}|\geq|V_{t+1}|= 1$ and $N^{+}_{D^{*}}(v_{1})=V_{2}\cup V_{3}\cup \ldots\cup V_{t}$, $N_{D^{*}}^{+}(v_{2})\supseteq V_{3}\cup \ldots \cup V _{t}\cup V_{t+1}$. So $v_{2}\rightarrow v_{t+2}$. Note that $v_{1}\rightarrow v_{2}$. Furthermore, $v_{t+2}\notin V_{1}$ since $D$ has no 2-cycle.

Note $v_{2}\rightarrow v_{t+2}$ and $v_{1}\nrightarrow v_{t+2}$. Then  $N^{-}_{D^{*}}(v_{t+2})=V_{2}\cup V_{3}\cup \ldots \cup V_{t}\cup V_{t+1}$.  Since $|V_{t+1}|=1$, we have $$d^{-}_{D^{*}}(v_{t+2})=|V_{2}\cup V_{3}\cup \ldots \cup V_{t}\cup V_{t+1}|=d_{D^{*}}^{+}(v_{1})+1=\delta^{+}(D^{*})+1\mbox{.}$$
Let $N_{D^{*}}^{+}(v_{t+2})=V_{t+3}\cup\ldots\cup V_{g}$. By Lemma \ref{121}, we have $$d^{++}_{D^{*}}(v_{t+2})\geq d^{+}_{D^{*}}(v_{g})\mbox{.}$$

{\noindent\bf Case 1.}  $v_{g}\neq v_{1}$.

We see that  $d^{++}_{D^{*}}(v_{t+2})\geq d^{+}_{D^{*}}(v_{g})\geq \delta^{+}(D^{*})+1\geq d_{D^{*}}^{-}(v_{t+2})$, and $v_{t+2}$ is a Sullivan-1 vertex of $D^{*}$.

Note $v_{t+2}\notin V_{1}$. Then $d_{D^{*}}^{+}(v_{t+2})\geq \delta^{+}(D^{*})+1\geq d_{D^{*}}^{-}(v_{t+2})$.
Thus $d_{D^{*}}^{++}(v_{t+2})+d_{D^{*}}^{+}(v_{t+2})\geq2 d_{D^{*}}^{-}(v_{t+2})$ and $v_{t+2}$ is a Sullivan-2 vertex of $D^{*}$.

Clearly, all the vertices of $V_{t+2}$ are Sullivan-$i$ vertices of $D^{*}$ for $i\in\{1,2\}$. By Corollary \ref{1.5}, the tournament $S_{t+2}$ has a Sullivan-1  vertex and a Sullivan-2 vertex, say $v_{t+2}$ and $v'_{t+2}$, respectively. By Lemma \ref{3.1}, $v_{t+2}$
is another Sullivan-1  vertex of $D$ and $v'_{t+2}$ is another Sullivan-2 vertex of $D$.

{\noindent\bf Case 2.}  $v_{g}=v_{1}$.

Now $N_{D^{*}}^{+}(v_{t+2})=V_{t+3}\cup\ldots\cup V_{r}\cup V_{1}$.

First we show that $v_{2}$ is another Sullivan-1 vertex of $D^{*}$.
Recall $v_{2}\rightarrow v_{t+2}$. Then $V_{2}\Rightarrow V_{3}\cup V_{4}\cup\ldots\cup V_{t+2}$. Also $V_{t+2}\Rightarrow V_{t+3}\cup\ldots\cup V_{r}\cup V_{1}$. Then $v\in N^{+}(v_{2})\cup  N^{++}(v_{2})$ for any $v\in V(D^{*}-V_{2})$. Since $N_{D^*-V_{2}}^{+}(v_{2})\cap N_{D^*-V_{2}}^{-}(v_{2})=\emptyset$, we have  $N_{D^{*}-V_{2}}^{-}(v_{2})\subseteq N_{D^{*}-V_{2}}^{++}(v_{2})$. Then $d_{D^{*}-V_{2}}^{++}(v_{2})\geq d_{D^{*}-V_{2}}^{-}(v_{2})$. This means that all vertices of $V_{2}$ are Sullivan-1 vertices of $D^{*}$. By Corollary \ref{1.5}, the tournament $S_{2}$ has a Sullivan-1  vertex, say also $v_{2}$. By Lemma \ref{3.1}, then $v_{2}$
is another Sullivan-1 vertex of $D$.

To find another Sullivan-2 vertex of $D$, we consider the following two cases.

If $d_{D^{*}}^{+}(v_{t+2})\geq \delta^{+}(D^{*})+2$, we have $d^{++}_{D^{*}}(v_{t+2})\geq d^{+}_{D^{*}}(v_{g})=d^{+}_{D^{*}}(v_{1})= \delta^{+}(D^{*})$.
Thus $d_{D^{*}}^{++}(v_{t+2})+d_{D^{*}}^{+}(v_{t+2})\geq\delta^{+}(D^{*})+\delta^{+}(D^{*})+2\geq2(\delta^{+}(D^{*})+1)\geq2 d_{D^{*}}^{-}(v_{t+2})$ and hence all vertices of $V_{t+2}$ are Sullivan-2 vertices of $D^{*}$. By Corollary \ref{1.5}, the tournament $S_{t+2}$ has a Sullivan-2 vertex, say also $v_{t+2}$. By Lemma \ref{3.1}, $v_{t+2}$ is another Sullivan-2 vertex of $D$.

If $d_{D^{*}}^{+}(v_{t+2})\leq\delta^{+}(D^{*})+1$, we have $d_{D^{*}}^{+}(v_{t+2})=\delta^{+}(D^{*})+1$
since we note $v_{t+2}\notin V_{1}$. Also note $v_{2}\rightarrow v_{t+2}$. Then $N_{D^{*}}^{+}(v_{2})\supseteq V_{3}\cup\ldots\cup V_{t+2}$ and hence $N_{D^{*}}^{-}(v_{2})\subseteq V_{t+3}\cup\ldots\cup V_{r}\cup V_{1}$. So $$d_{D^{*}}^{-}(v_{2})\leq|V_{t+3}\cup\ldots\cup V_{r}\cup V_{1}|=d_{D^{*}}^{+}(v_{t+2})=\delta^{+}(D^{*})+1\mbox{.}$$
Let $N_{D^{*}}^{+}(v_{2})=V_{3} \cup\ldots\cup V_{k}$. By Lemma \ref{121}, we have $d^{++}_{D^{*}}(v_{2})\geq d^{+}_{D^{*}}(v_{k})$. Note $v_{2}\rightarrow v_{k}$ and $v_{1}\rightarrow v_{2}$. Then $v_{k}\notin V_{1}$ since $D$ has no 2-cycle. Note that
$$d^{++}_{D^{*}}(v_{2})\geq d^{+}_{D^{*}}(v_{k})\geq \delta^{+}(D^{*})+1,\,\, d^{+}_{D^{*}}(v_{2})\geq \delta^{+}(D^{*})+1\mbox{.}$$
Then $d_{D^{*}}^{++}(v_{2})+d_{D^{*}}^{+}(v_{2})\geq2(\delta^{+}(D^{*})+1)\geq2 d_{D^{*}}^{-}(v_{2})$ and hence all vertices of $V_{2}$ are  Sullivan-2 vertices of $D^{*}$. By Corollary \ref {1.5}, the tournament $S_{2}$ always has a Sullivan-2 vertex, say also $v_{2}$. By Lemma \ref{3.1}, $v_{2}$ is another Sullivan-2 vertex of $D$.\qed

\begin{corollary}\label{4.5}
$(a)$ Every round decomposable local tournament has a Sullivan-$i$ vertex for $i\in\{1,2\}$.

$(b)$ Every round decomposable local tournament with no vertex of in-degree zero has at least two Sullivan-i vertices for $i\in\{1,2\}$.
\end{corollary}

Two examples of round local tournament, which have exactly two Sullivan-$i$ vertices for $i\in\{1,2\}$, are illustrated in Figure 1.

\begin{figure}[h]
      \unitlength0.6cm
      \begin{center}
      \begin{picture}(8,6)

      \put(-3,1){\circle*{.3}}
      \put(-3,4){\circle*{.3}}
      \put(1,1){\circle*{.3}}
      \put(1,4){\circle*{.3}}

      \qbezier(-3,1)(-3,2.5)(-3,4)
      \qbezier(1,1)(1,2.5)(1,4)
      \qbezier(-3,4)(-1,4)(1,4)
      \qbezier(-3,1)(-1,1)(1,1)
      \qbezier(-3,4)(-1,2.5)(1,1)
      \qbezier(-3,1)(-1,2.5)(1,4)

      \put(-3,2.5){\vector(0,1){.8}}
      \put(1,2.5){\vector(0,-1){.8}}
      \put(-1,4){\vector(1,0){.8}}
      \put(-1,1){\vector(-1,0){.8}}
      \put(-1,2.5){\vector(4,3){.8}}
      \put(-1,2.5){\vector(4,-3){.8}}

      \put(-3.8,1){\small $v_1$}
      \put(-3.8,4){\small $v_2$}
      \put(1.33,1){\small $v_4$}
      \put(1.33,4){\small $v_3$}
      \put(-1.4,0){\small $D$}


     \put(9,5.5){\circle*{.3}}
     \put(9,5.8){\small $v'_3$}

     \put(7.5,1){\circle*{.3}}
     \put(7.5,0.5){\small $v'_1$}

     \put(6.5,3.5){\circle*{.3}}
     \put(5.8,3.5){\small $v'_2$}

     \put(10.5,1){\circle*{.3}}
     \put(10.6,0.5){\small $v'_5$}

     \put(11.5,3.5){\circle*{.3}}
     \put(11.7,3.5){\small $v'_4$}

     \qbezier(6.5,3.5)(7.75,4.5)(9,5.5)
     \put(7.75,4.5){\vector(4,3){.08}}

     \qbezier(7.5,1)(9,1)(10.5,1)
     \put(9,1){\vector(-1,0){.08}}

     \qbezier(10.5,1)(11,2.25)(11.5,3.5)
     \put(11,2.24){\vector(-1,-2){.08}}

     \qbezier(11.5,3.5)(10.25,4.5)(9,5.5)
     \put(10.25,4.5){\vector(1,-1){.08}}

     \qbezier(6.5,3.5)(7,2.25)(7.5,1)
     \put(7,2.25){\vector(-1,2){.08}}

     \qbezier(7.5,1)(8.25,3.25)(9,5.5)
     \put(8.25,3.25){\vector(1,3){.08}}

     \qbezier(6.5,3.5)(8.5,2.25)(10.5,1)
     \put(8.5,2.25){\vector(-2,1){.08}}

     \qbezier(7.5,1)(9.5,2.25)(11.5,3.5)
     \put(9.5,2.25){\vector(-2,-1){.08}}

     \qbezier(10.5,1)(9.75,3.25)(9,5.5)
     \put(9.74,3.25){\vector(-1,4){.08}}
     \put(8.7,0){\small $D'$}
\end{picture}
   \caption {$D$ is a strong round local  tournament which has exactly two Sullivan-2 vertices $v_{1}$, $v_{2}$.  $D'$ is a strong round local  tournament which has exactly two Sullivan-$i$ vertices $v'_{4}$, $v'_{5}$ for $i\in\{1,2\}.$}
   \end{center}
   \end{figure}
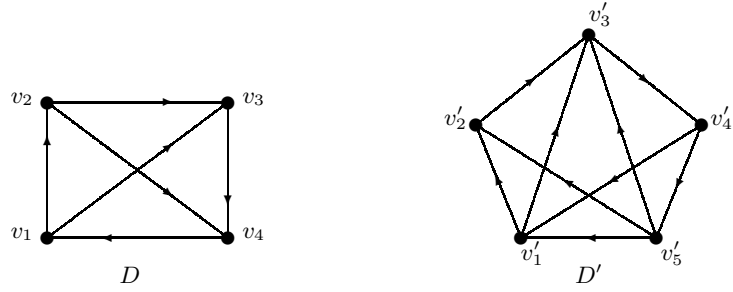

\section{In a non-round decomposable local tournament}

In this section, $D$ is always  a non-round decomposable local tournament, which is not a tournament. We also assume that $S$ is chosen with minimum cardinality among all minimal separating sets of $D$ satisfying that $D-S$ is not a tournament,
$D_{1},D_{2},\ldots,D_{p}$ is the acyclic ordering of the strong components of $D-S$, $D_{p+1},D_{p+2},\ldots,D_{p+q}$ is the acyclic ordering of the strong components of $D\langle S\rangle$, $D'_{1}, D'_{2}, D'_{3}$ is the semicomplete decomposition of $D-S$ and $D_{\lambda}$ is the initial component of $D'_{2}$. Clearly, in a local tournament $D$, the subdigraphs $S$, $D'_{1}$, $D'_{2}$ and $D'_{3}$ are all tournaments. See Figure 2. Let
$$A=N_{D'_{2}}^{+}(D_{1})=V(D_{\lambda})\cup V(D_{\lambda+1})\cup\ldots\cup V(D_{i}), \quad B=N_{S}^{+}(A)\mbox{,}$$
$$X=N^{+}_{ D-A}(A)=B\cup (V(D'_{2})-A)\cup V(D_{p})\mbox{,}$$
$$D'=D\langle S-B\rangle\mbox{,\,\,when $|X|\leq|S|$}.$$

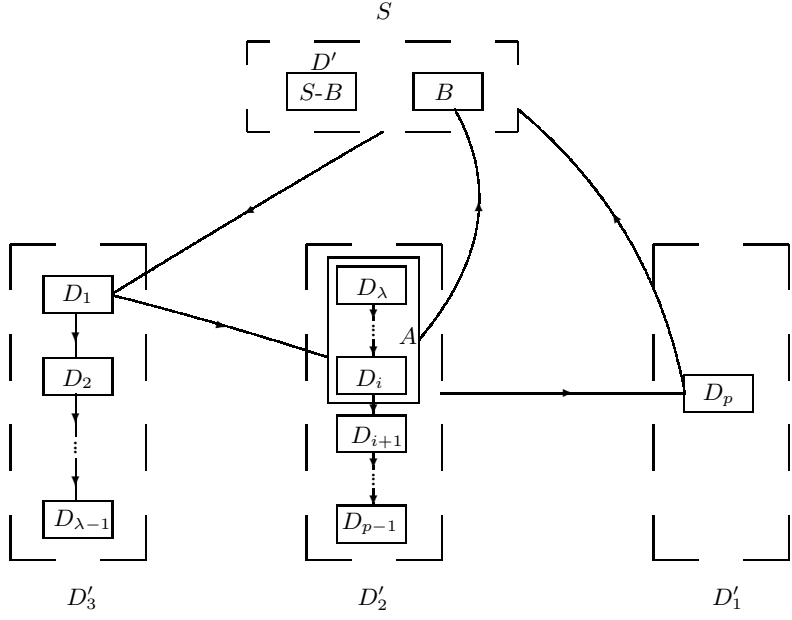
\begin{figure}[h]
      \unitlength0.6cm
   \begin{center}
      \begin{picture}(7,10)

      \put(2.85,12.5){\small$S$}
      \put(1.4,11.36){\small$D'$}
      \put(1.15,10.7){\small$S$-$B$}
      \put(4.15,10.7){\small$B$}

      \qbezier(4.6,10.5)(6,8)(3.82,5.4)
      \put(5.14,8.4){\vector(0,1){.08}}

      \put(-4.1,6.26){\small$D_{1}$}
      \put(-4.1,4.37){\small$D_{2}$}
      \put(-4.29,1.25){\small$D_{\lambda-1}$}
      \put(-4,-0.5){\small$D'_{3}$}

      \put(-3.8,5.4){\vector(0,-1){.08}}

      \put(3.37,5.3){\small$A$}
      \qbezier(-3,6.38)(-0.525,5.75)(1.79,5)
      \put(-0.525,5.7){\vector(4,-1){.08}}

      \put(-3.8,5){\line(0,1){1}}
      \put(-3.8,3.1){\circle*{.05}}
      \put(-3.8,3){\circle*{.05}}
      \put(-3.8,2.9){\circle*{.05}}

      \put(-3.8,1.8){\line(0,1){0.95}}
      \put(-3.8,2.2){\vector(0,-1){.08}}

      \put(-3.8,3.24){\line(0,1){0.95}}
      \put(-3.8,3.6){\vector(0,-1){.08}}

      \put(2.42,6.39){\small$D_{\lambda}$}
      \put(2.4,4.38){\small$D_{i}$}
      \put(2.3,3.1){\small$D_{i+1}$}
      \put(2.1,1.2){\small$D_{p-1}$}
      \put(2.45,-0.5){\small$D'_{2}$}

      \put(2.8,5.7){\circle*{.05}}
      \put(2.8,5.6){\circle*{.05}}
      \put(2.8,5.5){\circle*{.05}}

      \put(2.8,6){\line(0,1){0.2}}
      \put(2.8,5.9){\vector(0,-1){.08}}
      \put(2.8,5.2){\line(0,1){0.2}}
      \put(2.8,5.1){\vector(0,-1){.08}}

      \put(2.8,2.4){\circle*{.05}}
      \put(2.8,2.3){\circle*{.05}}
      \put(2.8,2.2){\circle*{.05}}

      \put(2.8,2.7){\line(0,1){0.2}}
      \put(2.8,2.6){\vector(0,-1){.08}}

      \put(2.8,1.88){\line(0,1){0.2}}
      \put(2.8,1.78){\vector(0,-1){.08}}

      \put(10.1,4.05){\small$D_{p}$}
      \put(10.3,-0.5){\small$D'_{1}$}

      \put(0,10){\dashbox(6,2)}
      \put(0.9,10.5){\framebox(1.5,0.8)}
      \put(3.7,10.5){\framebox(1.5,0.8)}

      \put(-5.24,0.5){\dashbox(3,7)}
      \put(-4.5,1){\framebox(1.5,0.8)}
      \put(-4.5,4.18){\framebox(1.5,0.8)}
      \put(-4.5,6){\framebox(1.5,0.8)}

      \put(1.32,0.5){\dashbox(3,7)}
      \put(2,0.9){\framebox(1.5,0.8)}
      \put(2,4.2){\framebox(1.5,0.8)}
      \put(2,6.2){\framebox(1.5,0.8)}
      \put(1.8,4){\framebox(2,3.2)}
      \put(2, 2.9){\framebox(1.5,0.8)}
      \put(2.8,3.7){\line(0,1){0.5}}

      \put(2.8,3.8){\vector(0,-1){.08}}

      \put(9,0.5){\dashbox(3,7)}

      \put(9.7,3.8){\framebox(1.5,0.8)}

      \qbezier(4.3,4.2)(7.15,4.2)(9.7,4.2)
      \put(7.15,4.2){\vector(1,0){.08}}

      \qbezier(6,10.5)(9,8)(9.7,4.2)
      \put(8.17,8.1){\vector(-2,3){.08}}

      \qbezier(3,10)(0,8.25)(-3,6.4)
      \put(0,8.23){\vector(-2,-1){.08}}
   \end{picture}
   \caption{The structure of a non-round decomposable local tournament $D$ in Section 4.}
   \end{center}
   \end{figure}

\begin{lemma}\label{555}
For any strong component $D_{j}$ of $D'_{3}$, $d_{D-D_{j}}^{+}(D_{j})\geq |S|$.
\end{lemma}
\noindent{\it Proof.} Let  $D_{j}\subseteq D'_{3}$ be arbitrary. Since $N_{D-D_{j}}^{+}(D_{j})$ is a separating set of $D$ and $D-N_{D-D_{j}}^{+}(D_{j})$ is not a tournament, we have $d_{D-D_{j}}^{+}(D_{j})\geq|S|$ by the choice of $S$.\qed

\vskip 3mm

\begin{lemma}\label{111} 
If $|X|<|S|$, then $S-B\Rightarrow A$ and $S-B\Rightarrow D'_{3}$.
\end{lemma}
\noindent{\it Proof.} Note that $X$ is also a separating set of $D$. By the choice of $S$, combining with $|X|<|S|$, we see that $D-X=D\langle V(D'_{3})\cup A \cup (S-B)\rangle$ is a tournament. So any vertex  of $S-B$ is adjacent to the vertices of $A$ and $D'_{3}$. Since $N_{S-B}^{+}(A)=\emptyset$, we have $S-B\Rightarrow A$. Also, an arc from $D'_{3}$ to $S$ implies that $D'_{3}$ and $D'_{1}$ are adjacent, which contradicts Theorem \ref{2.2} ($c$). So $S-B\Rightarrow D'_{3}$.\qed

\begin{lemma}\label{222} If $|X|<|S|$, then for any $v\in V(D')=S-B$,

 $(a)$ $X\subseteq N_{D-D'}^{+}(v)\cup N_{D-D'}^{++}(v)$.

 $(b)$ $N^{-}_{D-D'}(v)\subseteq N^{++}_{D-D'}(v)\subseteq X $.

 $(c)$ $d_{D-D'}^{+}(v)\geq d_{D-D'}^{-}(v)+2$.
\end{lemma}
\noindent{\it Proof.} $(a)$ Recall that  $X=B\cup (V(D'_{2})-A)\cup V(D_{p})$. Since $|X|<|S|$, we have $|B|\leq|X|-1\leq |S|-2$ and hence $|S-B|\geq2$. By Lemma \ref{111}, for any $v\in V(D')=S-B$, we have  $v\Rightarrow A \Rightarrow D'_{2}-A$ and  $v\Rightarrow A\Rightarrow D_{p}$. Since $B=N_{S}^{+}(A)$, for any $y\in B$, there exists a vertex $x\in A$ such that $x\rightarrow y$ and hence $v\rightarrow x\rightarrow y$. So any vertex  of $X$ either belongs to $N_{D-D'}^{+}(v)$ or belongs to $N_{D-D'}^{++}(v)$, i.e., $X\subseteq N_{D-D'}^{+}(v)\cup N_{D-D'}^{++}(v)$.

$(b)$ Recall that $D'=D\langle S-B\rangle$. By the definition of $A$, $B$ and $X$, we have $V(D-D')-X=A\cup V(D'_{3})$. For any $v\in V(D')$, by Lemma \ref{111}, we have $v\Rightarrow A$ and $v\Rightarrow D'_{3}$. Then $$N^{+}_{D-D'}(v)\supseteq A\cup V(D'_{3})=V(D-D')-X$$ and hence $N^{-}_{D-D'}(v)\subseteq X$ and $N^{++}_{D-D'}(v)\subseteq X.$
By $(a)$, $N^{-}_{D-D'}(v)\subseteq X\subseteq N_{D-D'}^{+}(v)\cup N_{D-D'}^{++}(v)$. Since $ N_{D-D'}^{+}(v)\cap  N_{D-D'}^{-}(v)=\emptyset$, we have  $N^{-}_{D-D'}(v)\subseteq N^{++}_{D-D'}(v)\subseteq X$.

$(c)$ By Lemma \ref{555}, we have $d^{+}_{D-D_{1}}(D_{1})\geq|S|$. Then $|A\cup V(D'_{3})|=|D_{1}|+d^{+}_{D-D_{1}}(D_{1})\geq1+|S|\geq|X|+2$ since $|X|<|S|$. Note that $A\cup V(D'_{3})\subseteq N^{+}_{D-D'}(v)$ and $N^{-}_{D-D'}(v)\subseteq X$. Then $$|X|+2\leq |A\cup V(D'_{3})|\leq d^{+}_{D-D'}(v),\,\,\,\,d^{-}_{D-D'}(v)\leq |X|.$$
So $d_{D-D'}^{+}(v)\geq d_{D-D'}^{-}(v)+2.$ \qed

\begin{lemma}\label{333} 
If $|X|<|S|$, then a Sullivan-i vertex of $D'$ is a Sullivan-i vertex of $D$ for $i\in\{1,2\}$.
\end{lemma}
\noindent{\it Proof.} By Lemma \ref{222} $(c)$,  for any $v\in V(D')$, $d^{+}_{D-D'}(v)\geq d^{-}_{D-D'}(v)+2$. Combining with  Lemma \ref{222} $(b)$, we have $d^{++}_{D-D'}(v)+d^{+}_{D-D'}(v)\geq2d^{-}_{D-D'}(v)+2$. By  Corollary \ref {1.5}, the tournament $D'$ always has a Sullivan-1 vertex and a Sullivan-2 vertex, say $v'$ and $v''$, respectively. Then $v'$ is a Sullivan-1 vertex of $D$ and $v''$ is a Sullivan-2 vertex of $D$.\qed

\begin{lemma}\label{444} 
If $|X|\geq|S|$, then a Sullivan-$i$ vertex  of  $D_{1}$ is a Sullivan-$i$ vertex of $D$ for $i\in\{1,2\}$.
\end{lemma}
\noindent{\it Proof.} Let $v$ be a Sullivan-$i$ vertex of $D_{1}$ for $i\in\{1,2\}$. We will prove that $v$ is a Sullivan-$i$ vertex of $D$ for $i\in\{1,2\}$.

Note that there is no arc from $D_{1}$ to $S$. Otherwise $D_{1}$ and $D'_{1}$ are adjacent which contradicts Theorem \ref{2.2} ($c$). So, $N_{D-D_{1}}^{+}(v)=N_{D-S-D_{1}}^{+}(v)$.  Combining with Theorem \ref{2.1} ($a$), we have
$$N_{D-D_{1}}^{+}(v)=N_{D-S-D_{1}}^{+}(v)=N_{D-S-D_{1}}^{+}(D_{1})=N_{D-D_{1}}^{+}(D_{1}).$$ Similarly, we have $N_{D-D_{1}}^{-}(v)=N_{D-D_{1}}^{-}(D_{1})$.
By Lemma \ref{555}, $$d^{+}_{D-D_{1}}(v)=d^{+}_{D-D_{1}}(D_{1})\geq|S|.$$ By the structure of  $D$  described in  Theorems \ref{2.5} and \ref{2.6}, we have  $N_{D-D_{1}}^{-}(D_{1})\subseteq S$ and $X\subseteq N_{D-D_{1}}^{++}(D_{1})=N_{D-D_{1}}^{++}(v)$. Now
$$d^{++}_{D-D_{1}}(v)\geq |X|\geq|S|,\,\,\,\, d^{-}_{D-D_{1}}(v)\leq|S|\mbox{.}$$

So $d^{++}_{D-D_{1}}(v)\geq|S|\geq d^{-}_{D-D_{1}}(v)$ and $d^{++}_{D-D_{1}}(v)+d^{+}_{D-D_{1}}(v)\geq2|S|\geq 2d^{-}_{D-D_{1}}(v)$. Since $v$ is also a Sullivan-$i$ vertex of $D_{1}$, we see that $v$ is a Sullivan-$i$ vertex of $D$.\qed

Now, we consider the existence of a Sullivan-1 vertex in $D$ and the number of Sullivan-1 vertices of $D$. In fact, the existence can be directly obtained from the following two results, which were proved by Wang, Yang and Wang \cite{Wang}, Li and Sheng \cite{Li quasi}, respectively.

\begin{lemma}\label{3.3}\cite{Wang} 
Let $D$ be a non-round decomposable locally semicomplete digraph. Then $D$ has a king.
\end{lemma}

\begin{proposition}\label{3.4} \cite{Li quasi}
 Let $D$ be an oriented graph. A king of $D$ is a Sullivan-1 vertex.
\end{proposition}

\begin{corollary}\label{3.5}  
Let $D$ be a non-round decomposable local tournament, which is not a tournament. Then $D$ has a Sullivan-1 vertex.
\end{corollary}

We consider primarily the number of  Sullivan-1 vertices.

\begin{theorem}\label{3.6} 
Let $D$ be a non-round decomposable local tournament,  which is not a tournament. Then $D$ has at least two Sullivan-1 vertices.
\end{theorem}

\noindent{\it Proof.} Recall that $A=N_{D'_{2}}^{+}(D_{1})=V(D_{\lambda})\cup V(D_{\lambda+1})\cup\ldots\cup V(D_{i}),\,\,B=N_{S}^{+}(A),\,\,X=N^{+}_{ D-A}(A)=B\cup (V(D'_{2})-A)\cup V(D_{p})$. The structure of $D$ is illustrated in Figure 2.

For the case when $|X|\geq|S|$, let $v$ be a Sullivan-1 vertex of $D_{1}$. By Lemma \ref{444}, $v$ is a Sullivan-1 vertex of $D$. By the proof of Lemma \ref{3.3} (See reference \cite{Wang}), there exists a king either belonging to $D_{2}'$ or belonging to $S$ in $D$, say $v'$. By Proposition \ref {3.4}, $v'$ is a Sullivan-1 vertex of $D$. Clearly, $v'\neq v$. Then $v'$ is another Sullivan-1 vertex of $D$.

For the case when $|X|<|S|$, note that $|B|\leq|S|-2$. Recall that  $D'=D\langle S-B\rangle$. Let $v$ be a Sullivan-1 vertex of $D'$. By Lemma \ref{333}, $v$ is Sullivan-1 vertex of $D$. Next we will find another Sullivan-1 vertex of $D$.

If $A=D'_{2}$, we will show that a king of $D_{p}$, say $v'$, is another Sullivan-1 vertex of $D$. By Theorem \ref{2.6} ($a$), we have $D_{p}\Rightarrow S$. Combining with Lemma \ref{111},  $D_{p}\Rightarrow D'\Rightarrow A=D'_{2}$ and $D_{p}\Rightarrow D'\Rightarrow D'_{3}$. Then $v'$ is a king of $D$. Clearly $v'\neq v$. By Proposition \ref{3.4}, $v'$ is another Sullivan-1 vertex of $D$.

Now $A\varsubsetneqq  D'_{2}$. Let $D''=D'-v$ and $u$ be a Sullivan-1 vertex of  $D''$. Assume $u\in V(D_{j})$, where $D_{j}$ is a strong  component of $S$. We consider the following two cases.

{\noindent\bf Case 1.} There exists no arc between $u$ and $D_{i+1}$. 

We will prove that $u$ is another Sullivan-1 vertex of $D$.

We claim that $N^{-}_{D-D'}(u)\varsubsetneqq  N^{++}_{D-D'}(u)$. By Lemma \ref{222} (b), we have $N^{-}_{D-D'}(u)\\\subseteq N^{++}_{D-D'}(u)$. We only need to prove that $N^{-}_{D-D'}(u)\neq N^{++}_{D-D'}(u)$. By Lemma \ref{111}, $u\Rightarrow A\Rightarrow D_{i+1}$. Combining with the fact that there exists no arc between $u$ and $D_{i+1}$, we have $D_{i+1}\subseteq N^{++}_{D-D'}(u)$ and
$D_{i+1}\nsubseteq N^{-}_{D-D'}(u)$. Then $N^{-}_{D-D'}(u)\neq N^{++}_{D-D'}(u)$.

Now $d^{++}_{D-D''}(u)\geq d^{++}_{D-D'}(u)\geq d^{-}_{D-D'}(u)+1\geq d^{-}_{D-D''}(u)$ since $N^{-}_{D-D'}(u)\varsubsetneqq  N^{++}_{D-D'}(u)$. Clearly $u\neq v$. Since $u$ is also a Sullivan-1 vertex of $D''$, we see that $u$ is another Sullivan-1 vertex of $D$.

{\noindent\bf Case 2.} There exists at least one arc between $u$ and $D_{i+1}$.

Let $g$ be a king of $D_{i+1}$. We will show that $g$ is another Sullivan-1 vertex.

We claim that $D_{i+1}\Rightarrow D_{j}$. Since $A$ is a separating set, we see that $D_{i+1},\ldots,D_{p},D_{p+1},\ldots,D_{p+q},D_{1},\ldots,D_{\lambda-1}$ is the acyclic ordering of the strong components of $D-A$. Since $D_{i+1}$ and $D_{j}$ are distinct strong components of $D-A$, by Theorem \ref{2.1} ($a$), we only need to prove that there exists at least one arc from $D_{i+1}$ to $D_{j}$. Now there exists at least one arc between $u$ and $D_{i+1}$ and hence there exists one arc between $D_j$ and $D_{i+1}$. However, by Theorem \ref{2.6} ($b$), an arc from $D_{j}$ to $D_{i+1}$ implies that $D_{1}\Rightarrow D_{i+1}$ which contradicts the definition of $A$. Then there exists no arc from $D_j$ to $D_{i+1}$ and hence there exists at least one arc from $D_{i+1}$ to $D_j$.

Note that $g\rightarrow u$ since $g\in V(D_{i+1})$ and $u\in V(D_{j})$. By Lemma \ref{111}, $u\Rightarrow A$ and $u\Rightarrow D'_{3}$. Hence $g\rightarrow u\Rightarrow A \,\,\,\mbox{and} \,\,\,g\rightarrow u\Rightarrow D'_{3}.$ By the structure of $D$ described in  Theorems \ref{2.5} and \ref{2.6}, we see that $g\Rightarrow D'_{2}-A-D_{i+1}$ and $g\Rightarrow D_{p}\Rightarrow S$. So $g$ is a king of $D$.  Clearly $g\neq v$. By Proposition \ref{3.4}, $g$ is another Sullivan-1 vertex of $D$.\qed

Next, we consider the existence of a Sullivan-2 vertex and the number of Sullivan-2 vertices in $D$.

\begin{lemma}\label{666} 
If $|X|=|S|$, then $D$ has at least two Sullivan-2 vertices of $D$.
\end{lemma}
\noindent{\it Proof.}    Let $v$ be a Sullivan-2 vertex of $D_{1}$. By Lemma \ref{444}, $v$ is a Sullivan-2 of $D$. To find another Sullivan-2 vertex of $D$, we consider the following two cases.

For the case when  $|X|=|S|$ and $D-X$ is not  a tournament, let $X$ be a minimal separating set of $D$ instead of $S$.
By Lemma \ref{333} and Lemma \ref{444}, there exists a Sullivan-2 vertex in $D$, say $u$. We can check $u\in S-B$ or $u\in X$ due to the new separating set $X$. Then $u\neq v$ since $v\in V(D_{1})$. Thus $u$ is another Sullivan-2 vertex of $D$.

For the case when $|X|=|S|$ and $D-X$ is  a tournament, note that $|B|\leq|X|-1=|S|-1$ since  $X=B\cup (V(D'_{2})-A)\cup V(D_{p})$. Recall that $D'=D\langle S-B\rangle$. Let  $u$ be a Sullivan-2 vertex of the tournament  $D'$.  We will show that $u$ is another Sullivan-2 vertex of $D$.

We claim $D'\Rightarrow A$ and $D'\Rightarrow D'_{3}$. Since $D-X=D\langle V(D'_{3})\cup A \cup
V(D')\rangle$ is a tournament,  we see that any vertex of $D'$ is adjacent to the vertices of $A$ and $D'_{3}$. Since $N_{D'}^{+}(A)=\emptyset$, we have $D'\Rightarrow A$. Also, an arc from $D'_{3}$ to $S$ implies that $D'_{3}$ and $D'_{1}$ are adjacent, which contradicts Theorem \ref{2.2} ($c$). So $D'\Rightarrow D'_{3}$.

By the definition of $A$, $B$ and $X$, we have $V(D-D')-X=A\cup V(D'_{3})$. Combining with $u\Rightarrow A$ and $u\Rightarrow D'_{3}$, we see that  $$ N^{+}_{D-D'}(u)\supseteq A\cup V(D'_{3})=V(D-D')-X$$ and hence $N^{-}_{D-D'}(u)\subseteq X,\,\,\,\,N^{++}_{D-D'}(u)\subseteq X.$

Also any vertex of $X=B\cup (V(D'_{2})-A)\cup V(D_{p})$ either belongs to $N_{D-D'}^{+}(u)$ or belongs to $N_{D-D'}^{++}(u)$ since $u\Rightarrow A \Rightarrow D_{2}^{'}-A$, $u\Rightarrow A\Rightarrow D_{p}$ and for any $g\in B$, there exists a vertex $h\in A$ such that $h\rightarrow g$. So
$X\subseteq N_{D-D'}^{+}(u)\cup N_{D-D'}^{++}(u)$. Since $ N_{D-D'}^{+}(u)\cap  N_{D-D'}^{-}(u)=\emptyset$, we have $N^{-}_{D-D'}(u)\subseteq N^{++}_{D-D'}(u)\subseteq X$ and hence $$d^{-}_{D-D'}(u)\leq d^{++}_{D-D'}(u)\leq |X|.$$
By Lemma \ref{555}, $d^{+}_{D-D_{1}}(D_{1})\geq|S|$.  Then $|A\cup V(D'_{3})|\geq|D_{1}|+d^{+}_{D-D_{1}}(D_{1})\geq1+|S|.$ Combining with  $  A\cup V(D'_{3})\subseteq N^{+}_{D-D'}(u)$, we have
$$d_{D-D'}^{+}(u)\geq |A\cup V(D'_{3})|\geq 1+|S|\geq1+|X|.$$

So $d^{++}_{D-D'}(u)+d^{+}_{D-D'}(u)\geq 2d^{-}_{D-D'}(u)+1$.
Clearly, $u\neq v$ since $ u\in S-B$ and  $v\in V(D_{1})$. The fact that   $u$ is a Sullivan-2 vertex of $D'$ implies that $u$ is another Sullivan-2 vertex of $D$.\qed

\begin{theorem}\label{3.7} Let $D$ be a non-round decomposable local tournament,  which is not a tournament. Then either $D$ has at least  two Sullivan-2 vertices or $D$ has a Sullivan-2 vertex $v$ satisfying  $d^{++}(v)+d^{+}(v)\geq 2d^{-}(v)+2$.
\end{theorem}
\noindent{\it Proof.} Recall that $A=N_{D'_{2}}^{+}(D_{1})=V(D_{\lambda})\cup V(D_{\lambda+1})\cup\ldots\cup V(D_{i})$, $B=N_{S}^{+}(A)$ and $X=B\cup (V(D'_{2})-A)\cup V(D_{p})$. We consider the following cases.

{\noindent\bf Case 1.} $|X|<|S|$.

Note that $|B|\leq |S|-2$. Recall that $D'=D\langle S-B\rangle$. Let $v$ be a Sullivan-2 vertex of $D'$. By Lemma \ref{333}, $v$ is a Sullivan-2 vertex of $D$. Let $D''= D'-v$ and $u$  be a Sullivan-2 vertex of the tournament $D''$. We will prove that $u$ is another Sullivan-2 vertex of $D$.

By Lemma \ref{222} $(b)$ and $(c)$, we see that  $d^{++}_{D-D'}(u)+d^{+}_{D-D'}(u)\geq2d^{-}_{D-D'}(u)+2$ and hence $d^{++}_{D-D''}(u)+d^{+}_{D-D''}(u)\geq 2d^{-}_{D-D''}(u)$. Since $u$ is also a Sullivan-2 vertex of $D''$, we have $d^{++}_{D''}(u)+d^{+}_{D''}(u)\geq 2d^{-}_{D''}(u)$.   Clearly $u\neq v$. Thus $d^{++}_{D}(u)+d^{+}_{D}(u)\geq 2d^{-}_{D}(u)$  and $u$ is another Sullivan-2 vertex of $D$.

{\noindent\bf Case 2.} $|X|=|S|$.

By Lemma \ref{666}, $D$ has at least two Sullivan-2 vertices.

{\noindent\bf Case 3.} $|X|>|S|$ and $|D_{1}|\geq2$.

By Lemma \ref{444}, a Sullivan-2 vertex of $D_{1}$ is a Sullivan-2 vertex of $D$. Since $D_{1}$ is strong and $|D_{1}|\geq2$, we see that $D_{1}$ has at least three vertices. By Corollary \ref{1.6}, $D_{1}$ has at least two Sullivan-2 vertices. These two  vertices are Sullivan-2 vertices of $D$.

{\noindent\bf Case 4.} $|X|>|S|$, $|D_{1}|=1$ and $D_{2}\subseteq D'_{2}$.

By Lemma \ref{444}, the only vertex $v$ of $D_{1}$ is a Sullivan-2 vertex of $D$. Note that $A=N^{+}_{D-D_{1}}(D_{1})$. By Lemma \ref{555}, $d_{D-D_{1}}^{+}(D_{1})=|A|\geq|S|$. We consider the following two cases. In the first case, we can find two Sullivan-2 vertices in $D$. In the second case, we can find a vertex $v$ satisfying  $d^{++}(v)+d^{+}(v)\geq 2d^{-}(v)+2$ in $D$.

If $|A|=|S|$, let $A$ be a minimal separating set of $D$ instead of $S$. By Lemma \ref{333} and Lemma \ref{444}, there exists a Sullivan-2 vertex in $D$, say $u$. We can check $u\in A$ or $u\in V(D_{i+1})$ due to the new separating set $A$. Then $u\neq v$ since $v\in V(D_{1})$. Thus $u$ is another Sullivan-2 vertex of $D$.

If $|A|>|S|$, we see that $N_{D-D_{1}}^{+}(v)=A$, $X\subseteq N_{D-D_{1}}^{++}(v)$ and $N_{D-D_{1}}^{-}(v)\subseteq S$. Then $d_{D-D_{1}}^{+}(v)=|A|\geq|S|+1$, $d_{D-D_{1}}^{++}(v)\geq |X|\geq|S|+1$ and $d_{D-D_{1}}^{-}(v)\leq |S|$ and hence $d_{D-D_{1}}^{++}(v)+d_{D-D_{1}}^{+}(v)\geq 2d_{D-D_{1}}^{-}(v)+2$. So $v$ is the desired vertex.

{\noindent\bf Case 5.} $|X|>|S|$, $|D_{1}|=1$ and $D_{2}\subseteq D'_{3}$.

By Lemma \ref{444}, the only vertex $v$ of $D_{1}$ is a Sullivan-2 vertex of $D$.
Let $$A'=N_{D'_{2}}^{+}(D_{2})=V(D_{\lambda})\cup V(D_{\lambda+1})\cup\ldots\cup V(D_{i'}), \quad B'=N_{S}^{+}(A')\mbox{,}$$
$$X'=N^{+}_{ D-A'}(A')=B'\cup (V(D'_{2})-A')\cup V(D_{p})\mbox{.}$$ The structure of $D$ is illustrated in Figure 3.
By Lemma \ref{555},  we have $$d_{D-D_{2}}^{+}(D_{2})\geq|S|.$$ Let $C=N_{D-D_{2}}^{+}(D_{2})=(D_{3}\cup\ldots\cup D_{\lambda-1})\cup A'$. Then $|C|\geq|S|.$

\begin{figure}[h]
      \unitlength0.6cm
      \begin{center}
      \begin{picture}(7,13)

       \put(2.85,12.5){\small$S$}
       \put(1.37,11.4){\small$D''$}
       \put(1.15,10.7){\small$S$-$B'$}
       \put(4.15,10.7){\small$B'$}

       \qbezier(4.6,10.5)(6,8)(3.82,5.4)
       \put(5.14,8.4){\vector(0,1){.08}}

       \put(-4.1,6.26){\small$D_{1}$}
       \put(-4.1,4.37){\small$D_{2}$}
       \put(-4.29,1.25){\small$D_{\lambda-1}$}
       \put(-4,-0.5){\small$D'_{3}$}

       \put(-3.8,5.4){\vector(0,-1){.08}}

      \put(3.26,5.3){\small$A'$}
      \qbezier(-3,4.56)(-0.525,4.58)(1.79,4.6)
      \put(-0.525,4.58){\vector(1,0){.08}}

      \put(-3.8,5){\line(0,1){1}}
      \put(-3.8,3.1){\circle*{.05}}
      \put(-3.8,3){\circle*{.05}}
      \put(-3.8,2.9){\circle*{.05}}

      \put(-3.8,1.8){\line(0,1){0.95}}
      \put(-3.8,2.2){\vector(0,-1){.08}}

      \put(-3.8,3.24){\line(0,1){0.95}}
      \put(-3.8,3.6){\vector(0,-1){.08}}

      \put(2.42,6.39){\small$D_{\lambda}$}
      \put(2.4,4.38){\small$D_{i'}$}
      \put(2.15,3.1){\small$D_{i'+1}$}
      \put(2.1,1.2){\small$D_{p-1}$}
      \put(2.45,-0.5){\small$D'_{2}$}

      \put(2.8,5.7){\circle*{.05}}
      \put(2.8,5.6){\circle*{.05}}
      \put(2.8,5.5){\circle*{.05}}

      \put(2.8,6){\line(0,1){0.2}}
      \put(2.8,5.9){\vector(0,-1){.08}}
      \put(2.8,5.2){\line(0,1){0.2}}
      \put(2.8,5.1){\vector(0,-1){.08}}

      \put(2.8,2.4){\circle*{.05}}
      \put(2.8,2.3){\circle*{.05}}
      \put(2.8,2.2){\circle*{.05}}

      \put(2.8,2.7){\line(0,1){0.2}}
      \put(2.8,2.6){\vector(0,-1){.08}}

      \put(2.8,1.88){\line(0,1){0.2}}
      \put(2.8,1.78){\vector(0,-1){.08}}

      \put(10.1,4.05){\small$D_{p}$}
      \put(10.3,-0.5){\small$D'_{1}$}

      \put(0,10){\dashbox(6,2)}
      \put(0.9,10.5){\framebox(1.5,0.8)}
      \put(3.7,10.5){\framebox(1.5,0.8)}

      \put(-5.24,0.5){\dashbox(3,7)}
      \put(-4.5,1){\framebox(1.5,0.8)}
      \put(-4.5,4.18){\framebox(1.5,0.8)}
      \put(-4.5,6){\framebox(1.5,0.8)}

      \put(1.32,0.5){\dashbox(3,7)}
      \put(2,0.9){\framebox(1.5,0.8)}
      \put(2,4.2){\framebox(1.5,0.8)}
      \put(2,6.2){\framebox(1.5,0.8)}
      \put(1.8,4){\framebox(2,3.2)}
      \put(2, 2.9){\framebox(1.5,0.8)}
      \put(2.8,3.7){\line(0,1){0.5}}
      \put(2.8,3.8){\vector(0,-1){.08}}

      \put(9,0.5){\dashbox(3,7)}

      \put(9.7,3.8){\framebox(1.5,0.8)}

      \qbezier(4.3,4.2)(7.15,4.2)(9.7,4.2)
      \put(7.15,4.2){\vector(1,0){.08}}

      \qbezier(6,10.5)(9,8)(9.7,4.2)
      \put(8.17,8.1){\vector(-2,3){.08}}

      \qbezier(3,10)(0,8.25)(-3,6.4)
      \put(0,8.23){\vector(-2,-1){.08}}

     \end{picture}
   \caption{The structure of a non-round decomposable local tournament $D$ in Case 5 of the proof of  Theorem \ref{3.7}.}
   \end{center}
   \end{figure}
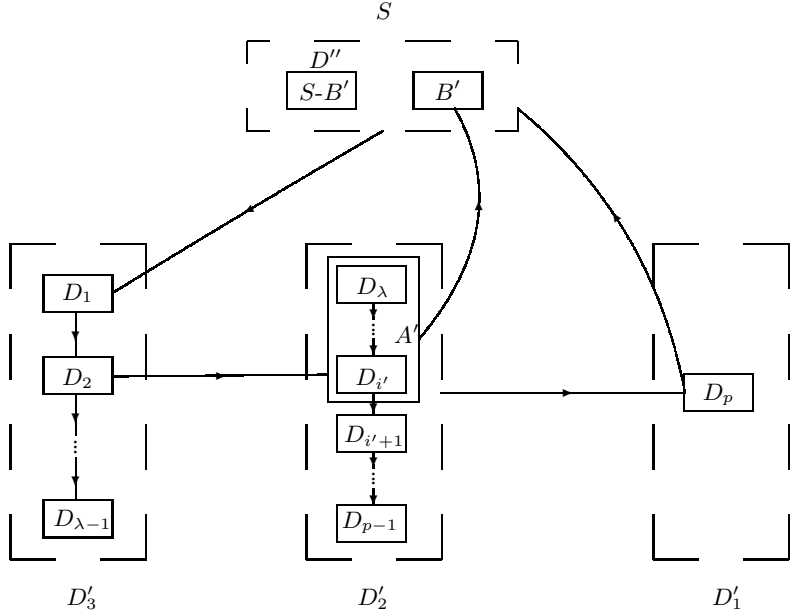

If $|X'|<|S|$, let $D''=D\langle S-B'\rangle$ and $u$ be a Sullivan-2 vertex of $D''$. Note that $u\neq v$ since $v\in V(D_{1})$.   Similarly to  the proof of  Lemma \ref{333}, we can show that $u$  is another Sullivan-2 vertices of $D$.

If $|X'|=|S|$, similarly to the proof of  Lemma \ref{666}, we can prove that $D$ has  two Sullivan-2 vertices.

If $|C|=|S|$, let $C$ be a minimal separating set of $D$ instead of $S$. By Lemma \ref{333} and Lemma \ref{444}, there exists a Sullivan-2 vertex in $D$, say $u$. We can check $u\in C$ or $u\in V(D_{i'+1})$ due to the new separating set $C$. Then  $u\neq v$ since $v\in V(D_{1})$. Then $u$ is another Sullivan-2 vertex of $D$.

If $|C|>|S|$ and $|X'|>|S|$, let $u$ be a Sullivan-2 vertex of $D_{2}$.  We will show that $u$ is another Sullivan-2 vertex of $D$. By the structure of $D$  described in  Theorems \ref{2.5} and \ref{2.6}, we have $N_{D-D_{2}}^{-}(u)\subseteq S\cup D_{1}$, $N_{D-D_{2}}^{+}(u)=C$, $X'\subseteq N_{D-D_{2}}^{++}(u)$. Then $d_{D-D_{2}}^{-}(u)\leq|S|+|D_{1}|=|S|+1$,
$d_{D-D_{2}}^{+}(u)=|C|\geq|S|+1,d_{D-D_{2}}^{++}(u)\geq |X'|\geq|S|+1.$
So $d_{D-D_{2}}^{++}(u)+d_{D-D_{2}}^{+}(u)\geq2(|S|+1)\geq 2d_{D-D_{2}}^{-}(u)$. Clearly $u\neq v$. Since $u$ is also a Sullivan-2 vertex of $D_{2}$, we see that $u$ is another Sullivan-2 vertex of $D$.

In any case, we find either two Sullivan-2 vertices or a vertex $v$ satisfying $d^{++}(v)+d^{+}(v)\geq 2d^{-}(v)+2$ in $D$. The proof of the theorem is complete.\qed

\begin{corollary}\label{3.8}
$(a)$ Every non-round decomposable local tournament which is not a tournament has a Sullivan-$i$ vertex for $i\in\{1,2\}$.

$(b)$ Every non-round decomposable local tournament which is not a tournament has at least  two Sullivan-1 vertices.

$(c)$ Every  non-round decomposable local tournament which is not a tournament either has at least two Sullivan-2 vertices or has a Sullivan-2 vertex $v$ satisfying $d^{++}(v)+d^{+}(v)\geq2d^{-}(v)+2$.
\end{corollary}


\section{Conclusion}

According to a full classification of local tournaments in Corollary \ref{1.8}, Corollary \ref{1.5}, Theorems \ref{1.61}, Corollary \ref{4.5} and Corollary \ref{3.8} imply the following theorem.

\begin{theorem} Let $D$ be a local tournament. Then the following holds for $D$.

$(a)$ Every local tournament $D$ has a Sullivan-$i$ vertex for $i\in\{1,2\}$.

$(b)$ Every local tournament $D$ with no vertex of in-degree zero has at least two Sullivan-1 vertices.

$(c)$ Every local tournament $D$ with no vertex of in-degree zero either has at least two Sullivan-2 vertices or has a Sullivan-2 vertex $v$ satisfying  $d^{++}(v)+d^{+}(v)\geq 2d^{-}(v)+2$.

\end{theorem}

\vskip 2mm



\begin{thebibliography}{99}
\bibitem{Bang-Jensen} {J. Bang-Jensen and G. Gutin. {\it Digraphs: Theory, Algorithms and Applications, 2nd edn.}, Springer-Verlag, London, 2009.}

\bibitem{Dean} {N. Dean and B. J. Latka. Squaring the tournament-an open problem. {\it Congr. Numer.}, 109 (1995) 73-80.}

\bibitem{Fisher} {D. C. Fisher. Squaring a tournament: a proof of Dean's conjecture. {\it J. Graph Theory}, 23(1) (1996) 43-48.}

\bibitem{Havet} {F. Havet and S. Thomass\'{e}. Median orders of tournaments: a tool for the second neighbourhood problem and Sumner' conjecture. {\it J. Graph Theory}, 35(4) (2000) 244-256.}

\bibitem{Fidler2} {D. Fidler and R. Yuster. Remarks on the second neighbourhood problem. {\it J. Graph Theory}, 55(3) (2007) 208-220.}

\bibitem{Ghazal} {S. Ghazal. Seymour's second neighborhood conjecture for tournaments missing a generalized star. {\it J. Graph Theory}, 71(1) (2012) 89-94.}

\bibitem{Kaneko} {Y. Kaneko and S. C. Locke. The minimum degree approach for Paul Seymour's distance 2 conjecture. {\it Congr. Numer}., 148 (2001) 201-206.}

\bibitem{Cohn} {Z. Cohn, A. Godbole, E. Wright Harkness and Y. Zhang. The number of Seymour vertices in random tournaments and digraphs. {\it Graphs Combin.}, 32(5) (2016) 1805-1816.}

\bibitem{Gutin}{G. Gutin and R. Li. Seymour's second neighbourhood conjecture for quasi-transitive oriented graphs. arxiv.org/abs/1704.01389.}

\bibitem{Chen} {G. Chen, J. Shen and R. Yuster. Second neighbourhood via first neighbourhood in digraphs. {\it Ann. Combin.}, 7(1) (2003) 15-20.}

\bibitem{Sullivan} {B. Sullivan. A summary of results and problems related to the Caccetta-H\"{a}ggkvist conjecture. arxiv.org/abs/math/0605646.}

\bibitem{Li quasi} {R. Li and B. Sheng. The second neighbourhood for quasi-transitive oriented graphs. {\it Acta Math. Sinica, English Series}, doi.org/10.1007/s10114-018-7287-3.}

\bibitem{Li tournament} {R. Li and B. Sheng. The second neighbourhood for bipartite tournaments.  {\it Discuss. Math. Graph Theory}, doi.org/10.7151/dmgt.2018.}

\bibitem{Bj}{J. Bang-Jensen. Locally semicompiete digraph: A generalization of tournaments. {\it J. Graph Theory}, 14 (1990) 371-390.}

\bibitem{Guo locally} {J. Bang-Jensen, Y. Guo, G. Guin and L. Volkmann. A classification of locally semicomplete digraphs. {\it Discrete Math.}, 167 (1997) 101-114.}

\bibitem{Li seymour}{R. Li and J. Liang. Seymour's second neighbourhood conjecture for local tournaments. preprinted.}

\bibitem{Con} {Y. Guo and L. Volkmann. Connectivity properties of locally semicomplete didraphs. {\it J. Graph Theory}, 18 (1994) 269-280.}

\bibitem{Wang} {R. Wang,  A. Yang and S.  Wang. Kings in locally semicomplete digraphs. {\it J. Graph Theory}, 63(4) (2010) 279-287.}
\end{thebibliography}
\end{document}